\documentclass[12pt]{article}

\usepackage{indentfirst}
\usepackage{amsfonts}
\usepackage{amsmath}
\usepackage{amssymb}
\usepackage{amsbsy}
\usepackage{graphicx}
\usepackage{epstopdf}
\usepackage{epic}

\newtheorem{dfn}{Definition}
\newtheorem{theorem}[dfn]{Theorem}
\newtheorem{lemma}[dfn]{Lemma}

\newtheorem{corollary}[dfn]{Corollary}
\newenvironment{pf}{\noindent{\bf Proof.}}
{\enspace\vrule height5pt depth0pt width5pt}
\def\deg{{\rm deg}}

\title{Deploying robots with two sensors in $K_{1,6}$-free graphs}

\begin{document}
\centerline{{\large\bf DEPLOYING ROBOTS WITH TWO SENSORS}}
\smallskip
\centerline{{\large\bf IN $K_{1,6}$-FREE GRAPHS}}
\bigskip
\bigskip
\centerline{{\bf Waseem Abbas}}
\smallskip
\centerline{{\bf Magnus Egerstedt}}
\smallskip
\centerline{School of Electrical and Computer Engineering}
\centerline{Georgia Institute of Technology}
\centerline{Atlanta, Georgia 30332, USA}
\bigskip
\centerline{{\bf Chun-Hung Liu}}
\smallskip
\centerline{{\bf Robin Thomas}%
\footnote{Partially supported by NSF under Grant No.~DMS-1202640.}}
\smallskip
\centerline{{\bf Peter Whalen}}
\smallskip
\centerline{School of Mathematics}
\centerline{Georgia Institute of Technology}
\centerline{Atlanta, Georgia  30332-0160, USA}
\bigskip

\begin{abstract}
\noindent
Let $G$ be a graph of minimum degree at least two with 
no induced subgraph isomorphic to $K_{1,6}$.
We prove that if $G$ is not isomorphic to  one of eight exceptional graphs,
then it is possible to assign  two-element subsets of $\{1,2,3,4,5\}$ to
the vertices of $G$ in such a way that for every $i\in\{1,2,3,4,5\}$ and
every vertex $v\in V(G)$ the label $i$ is assigned to $v$ or one of its
neighbors.
It follows that $G$ has fractional domatic number at least $5/2$.
This is motivated by a problem in robotics  and generalizes a result of
Fujita, Yamashita and Kameda who proved that the same conclusion holds for
all $3$-regular graphs.
\end{abstract}

\section{Introduction}

\def\junk#1{}
\junk{
Our work is motivated by a problem in robotics.
Assume that a network of robots is exploring 
 a broad terrain such as a polar ice cap, the surface of Mars, or the bottom of the ocean.
The robots can communicate, but there are physical limitations to which pairs of robots
can communicate with each other. Usually it is assumed that robots that are within 
certain distance can communicate, but we will not be concerned with physical
locations of the robots and instead will assume that the ability to communicate is
described by a graph $G$ with vertex-set the set of robots, where two robots are
adjacent if they can directly communicate.
It is reasonable to expect the robots to move over time, but we will assume that the
graph $G$ does not change.
If the robots do change locations over time, then our graph $G$ could be interpreted
as some ``minimum communications guarantee" in the sense that robots that are adjacent
in $G$ are required to be within communication range throughout their movement,
while other pairs may occasionally gain the ability to communicate. 
}

\junk{
The robots carry several types of sensors, but the robots are small and each can carry
only a limited number of those sensors. Let $r\ge1$ be an integer. We will assume that
each robot can carry at most $r$ sensors.
There is an important requirement that each robot must have access to the data from
at least one sensor of each type. In other words, for every robot $R$ and every type 
of sensor, either $R$ or one of its neighbors in $G$ must carry a sensor of that type.
The question we study is what is the maximum number of types of sensors that
can be deployed.
Our main result states that for $r=2$ and under some reasonable restrictions on $G$
it is possible to accommodate five types of sensors.
}

The problem under consideration in this paper is motivated by a problem encountered both in the  multi-agent robotics and  mobile sensor networks domains.  Common to both of these two application areas is a collection of agents  that are 
equipped with sensors of various types, used for tasks such as environmental modeling, exploration of unknown terrains, surveillance of remote locations, and the establishment of sensor coverage for the purpose of event detection. Due to the scale of the multi-robot network, the agents have to act based on locally available information, and under various such distributed coordinated schemes, e.g., \cite{Bullo}, the robots interact and communicate with each other in order to gain the information needed to make informed decisions. These interactions, in turn, define an information exchange network that allows us to model the agents as vertices and information exchange channels as edges in a graph.
The inter-agent interactions moreover allow the agents to complement each others' resources and capabilities; thus enhancing the collective functionality of the system. As a result, the underlying network topology of multi-robot networks plays a crucial role in achieving the system level objectives within the network in a distributed manner. 

As an example, 
consider an application in which a group of robots is deployed at some remote location for the purpose of environmental monitoring. Each robot needs to obtain information about $s$ different sensing modalities (e.g., temperature, humidity, barometric pressure, and so on). However, owing to certain constraints such as power limitations and hardware footprints, an individual robot can have a maximum of $r<s$ sensors installed on it. As a result, the robots need to collect data concerning the remaining $s-r$ sensing modalities from neighboring robots through the information exchange network. In other words, for every robot $v$ and every type of sensor, either $v$ or one of its neighboring robots must carry a sensor of that type.

As already stated, the 
 multi-robot network can be modeled as a graph $G$, in which the vertex set represents robots, and the edges
 correspond to the interactions among robots. 
Typically, a robot may transmit data to other robots lying within a certain Euclidean distance, say $R$, away from it. Thus, an edge is formed between nodes $v$ and $u$ whenever $\|v-u\| \leq R$. 
This results in an \emph{$R$-disk proximity graph} model of the network, which is the typical model employed when studying multi-robot networks. As such, any graph class under consideration must be rich enough to capture this model for it to be relevant to robotics. In such a graph, a disk of radius $R$, which represents the transmission or interaction range of the node, is associated with every node $v$ that lies at the center of the disk. An edge exists between $v$ and all such nodes that lie within the disk of $u$.   
$R$-disk graphs are one of the most frequently used models for the analysis of the network topology related aspects of multi-robot systems, 
wireless sensor networks, and other ad-hoc networks (e.g., see \cite{Mesbahi}). $R$-disk graphs are geometric graphs as the existence of edges between vertices depends on the geometric configuration of vertices. However, the geometric property of such graphs can be translated into a graph theoretic one. In fact, it can be shown that $R$-disk graphs are indeed $K_{1,6}$-free, and this key observation motivates the study of $K_{1,6}$-free graphs in multi-agent robotics.
 
 In this paper, we study what is the maximum number of sensors that can be accommodated in a multi-robot network if each robot can have at most two types of sensors. Our main result states that under some mild conditions, it is possible to assign two distinct labels to each vertex in a $K_{1,6}$-free graph such that a set of five distinct labels always exist in the closed neighborhood of every vertex in $G$.

The same problem arises in various situations of locating facilities in a network.
Let us assume that every vertex of a graph can access only resources located at
neighboring vertices or at the vertex itself.
Now if some resource (such as a file, a printer or other service) must be accessible 
from every vertex of the graph, then copies of that resource need to be distributed
over the network to form a ``dominating set".
If every vertex of the graph has the capacity to accommodate at most $r$
distinct resources, then asking for the maximum number of resources that can be 
made available to every vertex of the graph leads to the same mathematical
question as the problem of the previous paragraph.
 
Let us be more precise now.
By a {\em graph} we mean a finite, simple, undirected graph; that is, loops and parallel edges
are not allowed.
For a vertex $v$ of a graph $G$, we denote the set of neighbors of $v$ by $N(v)$, and define 
$N[v]$, the {\em closed neighborhood of $v$}, to be $N(v) \cup \{v\}$.
Let $r\ge1$ be an integer.
Let $f$ be a function that maps the vertices of $G$ to $r$-element
subsets of some set $X$.
We define $R(f)$ to be the union of $f(v)$ over all vertices $v$ of $G$.
Following~\cite{fyk} we say that $f$  is an {\em $r$-configuration on $G$} if for every $x\in R(f)$ 
 and every vertex $v\in V(G)$ we have $x\in f(u)$ for some $u\in N[v]$.
We define $D_r(G)$ to be the maximum of $|R(f)|$ over all $r$-configurations on $G$.
Thus given a graph $G$ and integer $r\ge1$ the problems of the previous two paragraphs ask 
for the value of $D_r(G)$ .

The parameter $D_1(G)$ is known in the literature as the {\em domatic number} of $G$.
It was introduced by Cockayne and Hedetniemi~\cite{CocHed} and has since then been
the subject of a large number of publications.
Obviously $D_1(G)$ is at most the minimum degree of $G$ plus one, but testing whether
$D_1(G)\ge k$ is NP-complete for all $k\ge 3$. 
(Testing $D_1(G)\ge2$ is easy, because
$D_1(G)\ge2$ if and only if $G$ has no isolated vertex.)
A $(1+o(1))\ln n$-approximation algorithm for $D_1(G)$ was found by
Feige, Halld\'orsson, Kortsarz and Srinivasan~\cite{FeiHalKorSri}, who also showed
that their approximation factor  is essentially best possible. 

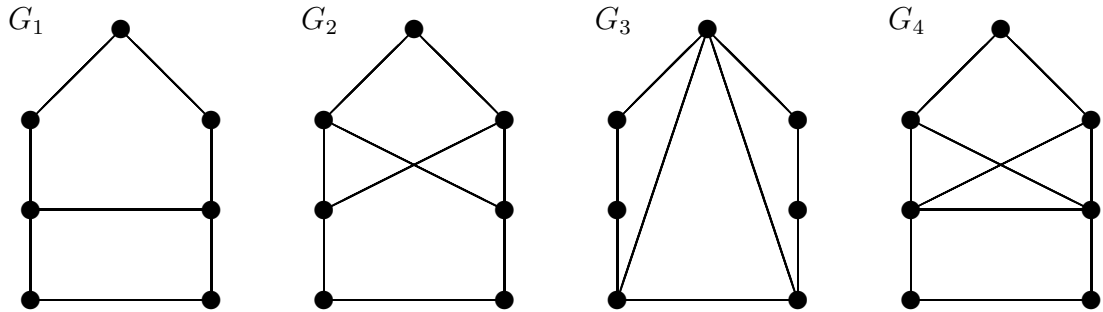
\begin{figure}
\unitlength=1.2mm
\begin{picture}(150, 30)(5,0)

\put(0,30){$G_1$}
\put(12.5,30){\circle*{2}}
\multiput(2.5,0)(0,10){3}{\circle*{2}}
\multiput(22.5,0)(0,10){3}{\circle*{2}}
\drawline[1000](12.5,30)(2.5,20)(2.5,0)(22.5,0)(22.5,20)(12.5,30)
\drawline[1000](2.5,10)(22.5,10)

\put(32.5,30){$G_2$}
\put(45,30){\circle*{2}}
\multiput(35,0)(0,10){3}{\circle*{2}}
\multiput(55,0)(0,10){3}{\circle*{2}}
\drawline[1000](45,30)(35,20)(35,0)(55,0)(55,20)(45,30)
\drawline[1000](35,20)(55,10)
\drawline[1000](35,10)(55,20)

\put(65,30){$G_3$}
\put(77.5,30){\circle*{2}}
\multiput(67.5,0)(0,10){3}{\circle*{2}}
\multiput(87.5,0)(0,10){3}{\circle*{2}}
\drawline[1000](77.5,30)(67.5,20)(67.5,0)(87.5,0)(87.5,20)(77.5,30)
\drawline[1000](67.5,0)(77.5,30)(87.5,0)

\put(97.5,30){$G_4$}
\put(110,30){\circle*{2}}
\multiput(100,0)(0,10){3}{\circle*{2}}
\multiput(120,0)(0,10){3}{\circle*{2}}
\drawline[1000](110,30)(100,20)(100,0)(120,0)(120,20)(110,30)
\drawline[1000](100,20)(120,10)
\drawline[1000](100,10)(120,20)
\drawline[1000](100,10)(120,10)
\end{picture}

\caption{Graphs $G_1$, $G_2$, $G_3$ and $G_4$.}    \label{Graph G_i}
\end{figure}

Fujita, Yamashita and Kameda proved in~\cite{fyk} that $D_2(G)\ge5$ for all $3$-regular graphs.
The purpose of this article is to generalize their result to a larger class of graphs, as follows.
We denote the cycle on $n$ vertices by $C_n$.
By {\em $C_4\cdot C_4$} we mean the graph obtained from two disjoint cycles on four
vertices by identifying a vertex in the first cycle with a vertex in the second cycle.
We denote by $G_1,G_2,G_3,G_4$ the graphs shown in Figure \ref{Graph G_i}.

\begin{theorem}
\label{thm:main}
Let $G$ be a graph of minimum degree at least two with no induced subgraph isomorphic to $K_{1,6}$.
If no component of $G$ is isomorphic to a member of $\{C_4,C_7,C_4\cdot C_4,K_{2,3},G_i: 1 \leq i \leq 4\}$, then $D_2(G)\ge5$.
\end{theorem}

As stated earlier, the generalization to $K_{1,6}$-free graphs is of interest in multi-agent robotics, because the class of $K_{1,6}$-free graphs includes the class of $R$-disk graphs.

For the sake of brevity let us define a {\em configuration} on a graph $G$ to mean a $2$-configuration
$f$ with $R(f)=\{1,2,3,4,5\}$.
Thus the conclusion of Theorem~\ref{thm:main} is equivalent to saying that $G$ has a configuration.
Our proof is algorithmic and gives a polynomial-time algorithm to find a configuration.
We say that a graph $G$ is {\em configurable} if it admits a configuration.
Theorem~\ref{thm:main} has the following two corollaries.

\begin{corollary}
If $G$ is a connected graph of minimum degree at least two with no induced subgraph isomorphic to $K_{1,6}$, and $G$ is not isomorphic to a member of $\{C_4,C_7,C_4\cdot C_4,K_{2,3},G_i: 1 \leq i \leq 4\}$, then for any positive integer $r$, $D_r(G)\ge \lfloor{5r}/{2}\rfloor$.
\end{corollary}

\begin{pf} 
Since $G$ has no isolated vertex, we have $D_1(G)\ge1$. Thus $G$ has a $1$-configuration $h$ with $R(h)=\{1,2\}$.
By Theorem~\ref{thm:main} the graph $G$ has a configuration, say $f$. 
For $v\in V(G)$ we define $g(v)$ to be the set of all pairs $(i,j)$,
where $i\in f(v)$ and $j\in\{1,2,\ldots,\lfloor r/2\rfloor\}$, and let $g'(v):=g(v)\cup h(v)$.
If $r$ is even, then $g$ is an $r$-configuration with $|R(g)|=5r/2$, and if
$r$ is odd, then $g'$ is an $r$-configuration with $|R(g')|=5(r-1)/2+2=\lfloor 5r/2\rfloor$, as desired.
\end{pf}
\bigskip

In the context of $R$-disk graphs, which are widely used to model inter-communication and information exchange among nodes in multi-robot and wireless sensor networks, we can restate the above result using the fact that $R$-disk graphs are always $K_{1,6}$-free, and can never be isomorphic to $K_{2,3}$, as shown in \cite{Leeuwen}.

\begin{corollary}
If $G$ is a connected $R$-disk graph of minimum degree at least two, and $G$ is not isomorphic to a member of $\{C_4,C_7,C_4\cdot C_4,K_{2,3},G_i: 1 \leq i \leq 4\}$, then for any positive integer $r$, $D_r(G)\ge \lfloor{5r}/{2}\rfloor$.
\end{corollary}

The {\em fractional domatic number} of a graph $G$, introduced in~\cite{Suo},
 is the supremum of $a/b$ such that $G$ has a $b$-configuration 
$f$ with $|R(f)|=a$. 
This is the optimum of the LP relaxation of the domatic number problem, and that justifies the name.
It follows that the supremum is attained.
Theorem~\ref{thm:main} implies that every graph that satisfies the hypotheses of the theorem has fractional domatic number at least $5/2$. 

\def\junk#1{}
\junk{
Given two positive integers $x$ and $y$, a function $f$ is called a {\it $(x,y)$-configuration function} for a graph $G$ if $f: V(G) \rightarrow [x]_y$ satisfies that $\cup_{u \in N[v]} f(u) = \{1,2, ..., x\}$ for all $v \in V(G)$, where $[x]_y$ is the collection of all $y$-subsets of $\{1,2,...,x\}$.  
Then the object of our problem is finding the maximum $x$ such that $G$ has a $(x,y)$-configuration function when $y$ is given.
It is clear that $G$ has a $(xt,yt)$-configuration function for every positive integer $t$ if $G$ has a $(x,y)$-configuration function.
A graph $G$ is called {\it $H$-free} if $G$ contains no induced subgraph which is isomorphic to $H$.
Fujita et al. \cite{fyk} proved that every $3$-regular graph has a $(5,2)$-configuration function.
In this paper, we prove a more general result that every $K_{1,6}$-free graph of minimum degree at least two has a $(5,2)$-configuration function, except four graphs.
Also, our constructive proof gives a polynomial time algorithm to find a $(5,2)$-configuration function for such graphs.
The interests of investigating $(5,2)$-configuration functions on $K_{1,6}$-free graphs is arisen from the well-known fact that unit disk graphs are $K_{1,6}$-free.
On the other hand, we give infinitely many examples that the condition of $K_{1,6}$-freeness can not be replaced by $K_{1,9}$-freeness.
In addition, we give infinitely many examples that the condition of forbidden graphs cannot be removed even for graphs of large minimum degree.
For the sake of simplicity, we omit the $(5,2)$ and refer to $f$ as just a \textit{configuration} function.
A graph $G$ is called {\it configurable} if there is a configuration function for $G$.
}

The paper is organized as follows.
In Section~\ref{sec:prelim} we prove some lemmas about extending a configuration from a subgraph of a graph.
In section 3 we prove the main theorem under the additional hypothesis that no two vertices of degree at least three are adjacent.
In section 4 we prove the main theorem and give two examples that show limitations to possible extensions.

\section{Preliminary lemmas}
\label{sec:prelim}
An $(\alpha, \beta)${\em-star} is the graph obtained by identifying one end of each of $\alpha$ paths of
length one and $\beta$ paths of length two.
In other words, the vertex-set may be labeled $\{w, x_i, y_j, z_j: 1 \leq i \leq \alpha, 1 \leq j \leq \beta\}$
so that the edge-set is $\{wx_i, wy_j, y_jz_j: 1 \leq i \leq \alpha, 1 \leq j \leq \beta\}$.
Note that an $(\alpha,0)$-star is isomorphic to $K_{1,\alpha}$.
We denote by $[5]^2$ the set of all two-element subsets of $\{1,2,3,4,5\}$.
If $G$ is a graph, $f: V(G) \rightarrow [5]^2$, and $v \in V(G)$, then we say that  \textit{$v$ is satisfied} 
with respect to $f$ if $\bigcup_{u \in N[v]} f(u) = \{1,2,3,4,5\}$.
When there is no danger of confusion we will omit the reference to $f$.


\begin{lemma} \label{P_2 helper}
Let $v_1v_2v_3v_4$ be a path of length three, and $f: \{v_1, v_4\} \rightarrow [5]^2$ with $f(v_1) \cap f(v_4)$ nonempty.
If $a,b \in \{1,2,3,4,5\} \setminus f(v_1)$, then $f$ can be extended to $\{v_1, v_2, v_3, v_4\}$ in such a way 
that $v_2$ and $v_3$ are satisfied and $f(v_2)=\{a,b\}$.
\end{lemma}

\begin{pf}
Without loss of generality, $f(v_1) = \{1, 2\}$, $1 \in f(v_4)$, and $f(v_2) = \{a, b\} = \{3, 4\}$.  Then setting $f(v_3) = \{2, 5\}$ completes the proof.
\end{pf}

\begin{lemma} \label{extendable}
Let $H$ and $S$ be disjoint subgraphs of a graph $G$, and let $\alpha,\beta\ge0$ be integers such that either
$\alpha + 3\beta \leq 9$ or $(\alpha, \beta) = (1,3)$.
Let $H$ be configurable and let $S$ be either a path of length at least two or an $(\alpha, \beta)$-star.
If every vertex of $S$ of degree one is adjacent to some vertex of $H$, 
then the subgraph of $G$ induced by $V(H) \cup V(S)$ is configurable. 
\end{lemma}

\begin{pf}
Let $f$ be a configuration  on $H$.
First, suppose that $S=v_1v_2...v_k$ is a path of length at least two (so $k \geq 3)$, and that the ends of $S$ are adjacent to vertices $x,y$ of $H$.
Note that $x$ and $y$ may be the same vertex.
There are three cases depending on the cardinality of $f(x)\cap f(y)$ and three cases depending on the residue of
$k$ modulo three.
Without loss of generality we may assume that $f(x)=f(y)=\{1,2\}$, or 
$f(x)=\{1,2\}$ and $f(y) = \{1,3\}$, or $f(x)=\{1,2\}$ and $f(y)=\{3,4\}$.
Then $f$ can be extended to $V(H) \cup V(S)$ according to the following table, where
 $t$ runs from $1$ through $\lfloor k/3 \rfloor -1$.
$$
\begin{tabular}{c | c | c | c | c | c | c | c}
\hline
$k$ (mod 3) & $f(x)$ & $f(v_{3t+1})$ & $f(v_{3t+2})$ & $f(v_{3t+3})$ & $f(v_{k-1})$ & $f(v_{k})$ & $f(y)$ \\ \hline
$0$ & $\{1,2\}$ & \{1,3\} & \{4,5\} & \{2,3\} & x & x & \{1,2\} \\ \hline
$0$ & $\{1,2\}$ & \{3,4\} & \{1,5\} & \{2,4\} & x & x & \{1,3\} \\ \hline
$0$ & $\{1,2\}$ & \{3,4\} & \{1,5\} & \{1,2\} & x & x & \{3,4\} \\ \hline
$1$ & $\{1,2\}$ & \{3,4\} & \{1,5\} & \{2,5\} & x & \{3,4\} & \{1,2\} \\ \hline
$1$ & $\{1,2\}$ & \{3,4\} & \{1,5\} & \{2,5\} & x & \{3,4\} & \{1,3\} \\ \hline
$1$ & $\{1,2\}$ & \{3,5\} & \{1,4\} & \{1,2\} & x & \{3,5\} & \{3,4\} \\ \hline
$2$ & $\{1,2\}$ & \{3,4\} & \{1,5\} & \{1,2\} & \{3,4\} & \{1,5\} & \{1,2\} \\ \hline
$2$ & $\{1,2\}$ & \{3,4\} & \{2,5\} & \{1,2\} & \{3,4\} & \{2,5\} & \{1,3\} \\ \hline
$2$ & $\{1,2\}$ & \{3,4\} & \{1,5\} & \{2,4\} & \{1,3\} & \{2,5\} & \{3,4\} \\ \hline
\end{tabular}
$$

Now we assume that $S$ is a $(\alpha, \beta)$-star, where $\alpha + \beta \geq 3$, $\alpha + 3 \beta \leq 9$, or $(\alpha, \beta) = (1,3)$.
Let $V(S) = \{w, x_i, y_j, z_j: 1 \leq i \leq \alpha, 1 \leq j \leq \beta\}$, $E(S) = \{wx_i, wy_j, y_jz_j: 1 \leq i \leq \alpha, 1 \leq j \leq \beta\}$, and $x_i$ is adjacent to $u_i$, where $u_i$ is in $H$, for all $1 \leq i \leq \alpha$, and $z_j$ is adjacent to $v_j$, where $v_j$ is in $H$, for all $1 \leq i \leq \beta$.

We say that $u_i$ {\em forbids} the set $f(u_i)$ and that $v_j$ {\em forbids} the three $2$-element subsets of $[5]-f(v_j)$.
We claim that there is an element of $[5]^2$ that is not forbidden by any $u_i$ or $v_j$. Indeed, this is clear if $\alpha + 3 \beta \leq 9$.
But if $\beta=3$, then the vertices $v_1,v_2,v_3$ collectively forbid at most eight sets, and hence the claim holds even when $\alpha=1$ and $\beta=3$.
We define $f(w)$ to be an element of $[5]^2$ that is not forbidden by any $u_i$ or $v_j$.
Furthermore, if $\beta=0$ and $|\bigcup_{i=1}^\alpha f(u_i)|\le3$, then we choose $f(w)$ disjoint from every $f(u_i)$.

If $\beta\ge1$, then we choose $f(x_i)$, $f(y_j)$ and $f(z_j)$ for $i=1,2,\ldots,\alpha$ and $j=1,2,\ldots,\beta-1$ in such a way that the
vertices $x_i,y_j,z_j$ are satisfied. 
Then $w$ sees at least three values under $f$ since any neighbor of $w$ already assigned a value does not have the exact same assignment as $w$.  So by Lemma \ref{P_2 helper} applied to the path $wy_\beta z_\beta v_\beta$ we can assign $f(y_\beta)$ and $f(z_\beta)$ in such a way that
 $y_\beta,z_\beta$ and $w$ are satisfied.
This completes the case $\beta\ge1$.

So we may assume $\beta = 0$. 
We assign $f(x_i)$ for $i=1,2,\ldots,\alpha$ such that $x_i$ is satisfied, $f(x_i)\cap f(w)=\emptyset$, and, if possible, not all $f(x_i)$ are the same.
Then $w$ is satisfied, unless the sets $f(x_i)$ are all equal, and so
from the symmetry we may assume that $f(w)=\{1,2\}$ and $f(x_i)=\{3,4\}$ for all $i=1,2,\ldots,\alpha$.
But then the choice of $f(x_i)$ implies that $f(u_i)\subseteq\{1,2,5\}$, contrary to the choice of $f(w)$.
\end{pf}

\begin{lemma} \label{degree two path}
Let $G$ be a graph, and let $P=xv_1v_2v_3y$ be a path in $G$.
If $x$ is adjacent to $y$, then let $H:=G\backslash \{v_1,v_2,v_3\}$; otherwise let $H$ be the graph obtained from
$G\backslash \{v_1,v_2,v_3\}$ by adding the edge $xy$.
If $H$ is configurable, then $G$ is configurable.
\end{lemma}

\begin{pf}
Let $f$ be a configuration on $H$.
We shall extend $f$ to $V(G)$.
If $f(x)=f(y)$, say $f(x)=\{1,2\}$, then $H \setminus xy$ is also configurable, so we can extend $f$ to $V(G)$ by Lemma \ref{extendable}.
So we may assume that $f(x) \neq f(y)$; that is, $\lvert f(x) \cup f(y) \rvert \geq 3$.
Define $g: V(G) \rightarrow [5]^2$ by $g(v_1)=f(y)$, $g(v_3)=f(x)$, 
let $g(v_2)$ be a $2$-element subset of $[5]^2$ containing $\{1,2,3,4,5\} \setminus (f(x) \cup f(y))$,
and let $g(v)=f(v)$ for all $v \in V(G) \setminus \{v_1, v_2, v_3\}$.
Then it is clear that $g$ is a configuration on $G$.
\end{pf}

\bigskip

Let $G$ be a graph and $v$ a vertex of $G$.
Let $f$ be a function mapping $V(G)$ to $[5]^2$ and $c \in [5]$.
Then we say that $v$ is {\it missing} $c$ if $c \not \in \bigcup_{u \in N[v]}f(u)$.

\begin{lemma} \label{extend a tailed cycle}
Let $H$ be $C_4$, $C_7$ or a configurable graph, and let $u_0$ be a vertex of $H$.
Let $G$ be a graph, where $V(G) = V(H) \cup \{u_i, w_j: 1 \leq i \leq k, 1 \leq j \leq m\}$ and $E(G) = E(H) \cup \{u_i u_{i+1}, u_kw_1, w_j w_{j+1}, w_mw_1: 0 \leq i \leq k-1, 1 \leq j \leq m-1\}$ for some nonnegative integer $k$ and integer $m$ with $m \geq 3$.
Then $G$ is configurable.
\end{lemma}

\begin{pf}
By Lemma \ref{degree two path} we may assume that $k=0,1$ or $2$.
Let $C$ be the cycle $w_1w_2 ... w_m w_1$.  Since $H$ is $C_4$, $C_7$ or a configurable graph, we may satisfy every vertex of $H$ except possibly $u_0$ and $u_0$ is missing at most two colors.  So we may assume $f(u_0) = \{1, 2\}$ and that $u_0$ is missing $3$ and $4$.  Similarly we may choose $f$
on $C$ in such a way that every vertex of $C$ except possibly $w_1$  is satisfied, and that $w_1$ is missing at most two colors.  

If $k = 0$ we choose $f$ on $C$ so that  $f(w_1) = \{3, 4\}$ and the colors missing at $w_1$ are $1$ and $2$.  
If $k = 1$, we choose $f$ on $C$ so that $f(w_1) = \{2, 5\}$ and the colors missing at $w_1$ are $3$ and $4$. We set $f(u_1) = \{3, 4\}$.  
Finally, if $k = 2$, we choose $f$ on $C$ so that  $f(w_1) = \{2, 3\}$ and the colors missing at $w_1$ are $1$ and $5$. 
We set $f(u_1) = \{3, 4\}$ and $f(u_2) = \{1, 5\}$.  
\end{pf}

\begin{lemma} \label{small extendable}
Let $H$ be a configurable graph, and let $f$ be a configuration on $G$.
If $G$ is obtained from $H$ by either
\begin{itemize}
	\item adding a vertex $v$ and two edges $vx$ and $vy$ to $H$, where $x,y$ are vertices of $H$ and $f(x) \neq f(y)$, or
	\item adding two vertices $u,v$ and three edges $xu, uv, vy$ to $H$, where $x,y$ are vertices of $H$ and $f(x) \cap f(y) \neq \emptyset$,
\end{itemize}
then $f$ can be extended to $G$.
\end{lemma}

\begin{pf}
This is easy to verify.
\end{pf}

\bigskip

A graph $G$ is said to be obtained from a graph $H$ by {\it attaching} a path $P$ if $G$ is obtained from the disjoint union of $H$ and $P$ by adding two edges $v_1x$ and $v_ky$, where $v_1$ and $v_k$ are the ends of $P$, and $x,y$ are vertices of $H$. 
A graph $G$ is said to be obtained from a graph $H$ by {\it adding} a path $P$ if $G$ is obtained from the disjoint union of $H$ and $P$ by identifying one end of $P$ and $x$ and identifying the other end of $P$ and $y$, where $x$ and $y$ are distinct vertices of $H$.

\begin{lemma} \label{5 6 cycle}
Let $C$ be a cycle of length of five or six.
If $G$ is obtained from $C$ by adding a path of length two or three between two nonadjacent vertices in $C$, then $G$ is configurable.
\end{lemma}

\begin{pf}
Let $C=v_1v_2...v_kv_1$, and $P$ be the path in $G \setminus C$ where the end of $P$ is adjacent to vertices $u,v$ of $C$ in $G$.
If $C$ is $C_5$, then we define a function $f: V(C) \rightarrow [5]^2$ by $f(v_i) = \{i, i+3\}$ for each $i=1,2,3,4,5$, where the addition is  modulo five.
If $C$ is $C_6$, then define $f(v_1)=\{1,3\}, f(v_2)=\{2,4\}, f(v_3)=\{1,5\}, f(v_4)=\{2,3\}, f(v_5)=\{1,4\}, f(v_6)=\{2,5\}$.
So $f(x) \neq f(y)$ for all distinct vertices $x,y$ in $C$, and $f(x) \cap f(y) \neq \emptyset$ for all nonadjacent two vertices $x,y$ in $C$.
Hence $f$ can be extended to $G$ by Lemma \ref{small extendable} since $P$ is a path on one or two vertices.
\end{pf}

\begin{lemma} \label{add a C_5 with two tails}
Let $x,y$ be vertices of a configurable graph $H$, let $C = v_1v_2...v_5v_1$ be a cycle of length five, and let $P= u_1 u_2 ... u_p$ and $Q = w_1 w_2 ... w_q$ be paths, where $p, q \in \{1,2\}$.
Assume that $H,C,P$ and $Q$ are pairwise disjoint.
If $G$ is the graph with $V(G) = V(H) \cup V(C) \cup V(P) \cup V(Q)$ and $E(G) = E(H) \cup E(C) \cup E(P) \cup E(Q) \cup \{xu_1, u_pv_1, yw_1, w_qv_3\}$, then $G$ is configurable.
\end{lemma}

\begin{pf}
Let $f$ be a configuration on $H$.
We shall extend $f$ to $G$.
If $f(x) \cap f(y)$ is nonempty, say $1 \in f(x) \cap f(y)$, then let $a,b$ are two distinct numbers in $\{1,2,3,4,5\} \setminus (f(x) \cup f(y))$, and define $f(v_1) = \{1,a\}$ and $f(v_3) = \{1,b\}$.
If $f(x)$ is disjoint from $f(y)$, say $f(x) = \{1,2\}$ and $f(y) = \{3,4\}$, then define $f(v_1) = \{1,3\}$ and $f(v_3) = \{1,4\}$.
Without loss of generality, we may assume that $a=3$ and $b=4$.
Then we further define $f(v_2)=\{2,5\}, f(v_4)=\{3,5\}$ and $f(v_5)=\{2,4\}$ so that every vertex of $C$ is satisfied.
By Lemma \ref{small extendable}, there is a way to define $f$ on $V(P) \cup V(Q)$ such that $f$ is a configuration  on $G$.
\end{pf}

\bigskip

Let us recall that the graph $C_4\cdot C_4$ was defined in the Introduction.

\begin{lemma} \label{add a path to cycle}
Let $G$ be a graph obtained by attaching a path $P=v_1v_2...v_k$ to a cycle $C$ with $v_1$ adjacent to $x$ and $v_k$ adjacent to $y$, for some vertices $x,y$ in $C$, where $k \geq 3$.
If $G$ is not isomorphic to $C_4\cdot C_4$ or $G_1$, then $G$ is configurable.
\end{lemma}

\begin{pf}
If $x$ is adjacent to $y$ in $C$, then $G$ is a cycle with a chord.
So $G$ is configurable when the cycle has length not four or seven.
It is easy to check that $G$ is configurable when the cycle has length four.
And since $G$ is not isomorphic to $G_1$, $G$ is also configurable when the cycle has length seven by Lemma \ref{5 6 cycle}.
So we may assume that $x$ is not adjacent to $y$ in $C$.
In other words, either $x$ equals $y$, or $x$ and $y$ are nonadjacent.

If the length of $C$ is not four or seven, then this lemma follows directly from Lemma \ref{extendable}.
So we may assume that the length of $C=u_1 u_2 ... u_{\lvert C \rvert} u_1$ is four or seven.
Also, we may assume that $3 \leq k \leq 5$ by Lemma \ref{degree two path}.
Without loss of generality, we assume that $x=u_1$.

{\bf Case 1:} $C = C_4$ and $x=y$.
Then $k=4$ or $5$ since $G$ is not isomorphic to $C_4\cdot C_4$.
So $G$ is isomorphic to the graph obtained by attaching a path of order three to  $C_5$ or  $C_6$, and hence $G$ is configurable by Lemma \ref{extendable}.

{\bf Case 2:} $C=C_4$ and $x \neq y$.
We may assume that $y = u_3$.
If $k=3$ or $5$, then $u_1 v_1 v_2 ... v_k u_3 u_2 u_1$ is a cycle of length six or eight, so it is configurable, and there is a configuration $f$ on it.
Then we can extend $f$ to $G$ by assigning that $f(u_3) = f(u_1)$, so $G$ is configurable.
If $k=4$, then we define a configuration  on $G$ by $f(u_1) = \{1,2\}, f(u_2) = \{3,5\}, f(u_3) = \{3,4\}, f(u_4) = \{2,5\}, f(v_1) = \{1,4\}, f(v_2) = \{3,5\}, f(v_3) = \{2,5\}, f(v_4) = \{1,4\}$.

{\bf Case 3:} $C=C_7$ and $x=y$.
We may assume that $x=y=u_1$.
If $k=4$ or $5$, then $G$ is isomorphic to the graph obtained by attaching a path of order six to  $C_5$ or  $C_6$, so $G$ is configurable by Lemma \ref{extendable}.
If $k=3$, then we can define a configuration on $G$ by $f(u_1) = \{1,2\}, f(u_2) = \{3,4\}, f(u_3) = \{1,5\}, f(u_4) = \{2,3\}, f(u_5) = \{1,4\}, f(u_6) = \{2,5\}, f(u_7) = \{3,4\}, f(v_1) = \{1,5\}, f(v_2) = \{3,4\}, f(v_3) = \{2,5\}$.

{\bf Case 4:} $C=C_7$, $x=u_1$ and $y=u_6$.
If $k=3$ or $5$, then $G$ is isomorphic to the graph obtained by attaching a path of order four to  $C_6$ or  $C_8$, so $G$ is configurable by Lemma \ref{extendable}.
If $k=4$, then we can define a configuration on $G$ by $f(u_1) = \{1,2\}, f(u_2) = \{3,4\}, f(u_3) = \{3,5\}, f(u_4) = \{1,2\}, f(u_5) = \{4,5\}, f(u_6) = \{3,4\}, f(u_7) = \{3,5\}, f(v_1) = \{1,5\}, f(v_2) = \{3,4\}, f(v_3) = \{2,5\}, f(v_4)=\{1,2\}$.

{\bf Case 5:} $C=C_7$, $x=u_1$ and $y=u_5$.
If $k=4$ or $5$, then $G$ is isomorphic to the graph obtained by attaching a path of order three to  $C_8$ or $C_9$, so $G$ is configurable by Lemma \ref{extendable}.
If $k=4$, then we can define a configuration on $G$ by $f(u_1) = \{1,2\}, f(u_2) = \{1,3\}, f(u_3) = \{4,5\}, f(u_4) = \{2,3\}, f(u_5) = \{1,2\}, f(u_6) = \{4,5\}, f(u_7) = \{3,4\}, f(v_1) = \{1,5\}, f(v_2) = \{3,4\}, f(v_3) = \{2,5\}$.
\end{pf}

\begin{lemma} \label{small case 2}
The graph $K_{2,4}$ is configurable.
\end{lemma}

\begin{pf}
Let $V(K_{2,4}) = \{x_1,x_2, y_1,y_2, y_3, y_4\}$, $E(K_{2,4}) = \{x_iy_j: 1 \leq i \leq 2, 1 \leq j \leq 4\}$.
We define a configuration  on $K_{2,4}$ by $f(x_1) = \{1,2\}, f(x_2) = \{3,4\}, f(y_1) = \{3,5\}, f(y_2) = \{4,5\}, f(y_3) = \{1,5\}, f(y_4) = \{2,5\}$.
\end{pf}

\begin{lemma} \label{add a path to small graph}
If a graph $G$ is obtained from $C_4\cdot C_4$ or $K_{2,3}$ by attaching a path, then $G$ is configurable.
\end{lemma}

\begin{pf}
First, we assume that $G$ obtained from  $C_4\cdot C_4$ by attaching a path $v_1v_2 ... v_k$, where $v_1$ is adjacent to $x$, $v_k$ is adjacent to $y$ for some vertices $x$, $y$ in $C_4\cdot C_4$.
We write the vertex set of $C_4\cdot C_4$ as $\{u_1, u_2, u_3, v, w_1, w_2,w_3\}$, where $vu_1u_2u_3v$ and $vw_1w_2w_3v$ are the two cycles in $C_4\cdot C_4$.

\noindent {\bf Case 1:} $x=y$.
By Lemma \ref{degree two path}, we may assume that $k = 2,3$ or $4$.
If $x=y=u_1$, then $G$ can be obtained from $C_3$ or $C_5$ by consecutively attaching a path of order three when $k=2$ or $4$, and $G$ has a spanning subgraph which is obtained from two disjoint $C_4$'s by attaching a path of order two when $k=4$, so $G$ is configurable by Lemma \ref{extendable} and Lemma \ref{extend a tailed cycle}.
Similarly, $G$ is configurable if both $x$ and $y$ are $u_3$, $w_1$ or $w_3$.
If $x=y=v_2$ and $k=2$ or $4$, then $G$ can be obtained from $C_3$ or $C_5$ by consecutively attaching a path of order three, so $G$ is configurable by Lemma \ref{extendable}.
If $x=y=u_2$ and $k=3$, then we define a configuration on $G$ as $f(v) = \{3,4\}, f(w_1) = \{1,3\}, f(w_2) = \{2,5\}, f(w_3) = \{1,4\}, f(u_1) = \{4,5\}, f(u_2) = \{1,2\}, f(u_3) = \{2,5\}, f(v_1) = \{1,3\}, f(v_2) = \{4,5\}, f(v_3) = \{2,3\}$.
Similarly, $G$ is configurable if $x=y=w_2$.
If $x=y=v$ and $k=2$ or $4$, then $G$ can be obtained from $C_3$ or $C_5$ by consecutively attaching a path of order three.
If $x=y=v$ and $k=3$, then we define a configuration by $f(v) = \{1,2\}, f(u_1) = \{1,3\}, f(u_2) = \{4,5\}, f(u_3) = \{2,3\}, f(v_1) = \{1,4\}, f(v_2) = \{3,5\}, f(v_3) = \{2,4\}, f(w_1) = \{1,5\}, f(w_2) = \{3,4\}, f(w_3) = \{2,5\}$.

\noindent {\bf Case 2:} $x \neq y$.
By Lemma \ref{degree two path}, we may assume that $k=0,1,2$.
When $k=0$, $G$ is obtained by adding an edge $xy$ to $C_4\cdot C_4$, and it is easy to show that $G$ is configurable.
When $k=1$, $x=v$, $y=u_2$, then define a configuration  on $G$ by $f(v) = \{1,2\}, f(u_1) = \{4,5\}, f(u_2) = \{3,4\}, f(u_3) = \{1,5\}, f(v_1) = \{2,5\}, f(w_1) = \{1,3\}, f(w_2) = \{4,5\}, f(w_3) = \{2,3\}$.
Similarly, $G$ is configurable if $k=1$, $x=w_1$ and $y=w_3$.
When $k=1$ and $x,y$ are not the case mentioned above, $G$ has a spanning subgraph which is $C_8$, or it can be obtained from either $C_5$ by attaching a path,  two disjoint $C_4$'s by adding an edge, or $C_5$ by attaching paths of order one or two, so $G$ is configurable by Lemma \ref{extendable}, Lemma \ref{extend a tailed cycle}, and Lemma \ref{small extendable}.

Now, we assume that $G$ obtained from $K_{2,3}$ by attaching a path $v_1v_2 ... v_k$, where $v_1$ is adjacent to $x$, $v_k$ is adjacent to $y$ for some vertices $x$, $y$ in $C_4\cdot C_4$.
We write $V(K_{2,3}) = \{u_1, u_2, w_1, w_2, w_3\}$ and $E(K_{2,3}) = \{u_iw_j: i=1,2, j =1,2,3\}$.

\noindent {\bf Case 3:} $x=y$.
By Lemma \ref{degree two path}, we may assume that $k=2,3,4$.
Then $G$ has a spanning subgraph which is obtained from either $C_3$ or $C_5$ by attaching a $(3,0)$-star, or $C_4\cdot C_4$ by attaching a path, or a cycle by attaching a $C_4$, so $G$ is configurable by Lemma \ref{extendable}, Lemma \ref{extend a tailed cycle}, Case 1 and Case 2.

\noindent {\bf Case 4:} $x \neq y$.
By Lemma \ref{degree two path}, we may assume that $k=0,1,2$.
If $x=u_1$, $y=u_2$ and $k=0$, then there is a configuration on $G$ defined by $f(u_1) = \{1,2\}, f(u_2) = \{3,4\}, f(w_1) = f(w_2)= f(w_3) = \{1,5\}$.
For other cases, $G$ contains a subgraph which is isomorphic to $K_{2,4}$ or $C_6$, or it can be obtained from either $C_3$ by attaching a path of order three,  $C_4\cdot C_4$ by adding an edge, $C_5$ or $C_6$ by attaching paths of order one or two, so $G$ is configurable by Lemma \ref{extendable}, Lemma \ref{small extendable}, Lemma \ref{small case 2}, Case 1 and Case 2.
\end{pf}

\section{A special case}

For a vertex $v$ of a graph $G$, we denote the degree of $v$ by $\deg_G(v)$.

\begin{lemma} \label{orientation}
For every graph $G$, there is an orientation of $E(G)$ such that each vertex $v$ has in-degree at least $\lfloor \deg_G(v)/2 \rfloor$.
\end{lemma}

\begin{pf}
We proceed by induction on $\lvert V(G) \rvert+\lvert E(G) \rvert$.
The lemma  obviously holds for the null graph. 
If $v$ is an isolated vertex of $G$, then the lemma follows by induction applied to $G\backslash v$.
If there is a vertex $v$ in $G$  of degree one, then, letting $u$ be the unique neighbor of $v$,  there is an orientation of $G \setminus uv$ such that the in-degree of each vertex $x$  is at least $\lfloor \deg_{G \setminus \{uv\}} (x)/2 \rfloor$ by the induction hypothesis, and then we can obtain a desired orientation of $G$ by orienting the edge $uv$  from $v$ to $u$.
So we may assume that $G$ has  minimum degree at least two, and hence $G$ contains a cycle $C=v_1v_2 ... v_k v_1$.
By the induction hypothesis, there is an orientation of $G \setminus E(C)$ such that the in-degree of each vertex $x$ is at least $\lfloor \deg_{G \setminus E(C)}/2 \rfloor$, and then we can obtain a desired orientation of $G$ by orienting the edges of $C$ to form a directed cycle.
This completes the proof.
\end{pf}

\bigskip

Note that the proof in Lemma \ref{orientation} gives a linear-time algorithm to find such an orientation.

\begin{lemma} \label{connectivity}
Let $H_1$ and $H_2$ be graphs, let $P$ be a path with at least one vertex, and let $v_1$ and $v_2$ be vertices of $H_1$ and $H_2$ respectively.  Let $G$ be the graph formed by taking the disjoint  union of $H_1$, $H_2$, and $P$ and identifying the first vertex of $P$ with $v_1$ and the last vertex of $P$ with $v_2$.  
Assume that $f_1$ and $f_2$ are functions mapping $V(H_1)$ and $V(H_2)$ to $[5]^2$, respectively, and that for $i=1,2$ the function $f_i$ satisfies every vertex 
of $H_i$ except possibly $v_i$.
If $|\bigcup_{u\in N(v_1)}{f_1(u)}| \geq 4$ and $|\bigcup_{u\in N(v_2)}{f_2(u)}| \geq 3$,
then $G$ is configurable.
\end{lemma}

\begin{pf}
Let $f'$ be the function defined to be $f_1$ on $H_1$ and $f_2$ on $H_2$.  Then $f'$ is a configuration for $G$ except possibly on $v_1$ and $v_2$ and $P$.  Suppose $|V(P)| \leq 2$.  Then we can permute the colors on $f_2$ so that $v_1$ and $v_2$ are satisfied, so we are done.  
If $|V(P)| = 3$, we may assume $f(v_1) = \{1, 2\}$ and $v_1$ is not missing a number except possibly $3$ and $f(v_2) =  \{4, 5\}$ and $v_2$ is not missing a number other than possibly $3$ and a number $c$.  
Then we set $f(u) = \{c, 3\}$ where $u$ is the middle vertex of $P$.  If $|V(P)| = 4$,  we apply Lemma \ref{P_2 helper}.  If $|V(P)| \geq 5$, we can reduce to one of the previous cases by applying Lemma \ref{degree two path}.
\end{pf}

\begin{lemma} \label{exceptional are almost satisfied}
Let $G$ be a graph and $v$ a vertex of $G$.
If $G$ is isomorphic to $C_4$, then there exists a function $f: V(G) \rightarrow [5]^2$ such that $v$ is satisfied and $\lvert \bigcup_{u \in N[v]} f(u) \rvert \geq 3$.
If $G$ is isomorphic to $C_7,C_4 \cdot C_4$ or $K_{2,3}$, then there exists a function $f: V(G) \rightarrow [5]^2$ such that $v$ is satisfied and $\lvert \bigcup_{u \in N[v]} f(u) \rvert \geq 4$.
\end{lemma}

\begin{pf}
This is easy to verify.
\end{pf}

\bigskip
We are now ready to prove an important special case of Theorem~\ref{thm:main}.

\begin{lemma} \label{min degree 2 max degree 5}
Let $G$ be a connected graph of maximum degree at most five and of minimum degree at least two with no two vertices of degree at least three adjacent.
If $G$ is not $C_4$, $C_7$, $C_4\cdot C_4$ or $K_{2,3}$, then $G$ is configurable.
\end{lemma}

\begin{pf}
Let $n$ be the order of $G$.
Suppose that $G$ is a minimum counterexample; that is, $G$ is not configurable, but $H$ is configurable for every graph $H$ with $\lvert V(H) \rvert + \lvert E(H) \rvert < \lvert V(G) \rvert + \lvert E(G) \rvert$ that satisfies the conditions of the lemma.

We note first that we may assume $G$ is 2-connected.
Otherwise we apply Lemma \ref{connectivity}, noting that each of the forbidden graphs except $C_4$ has the property that for every vertex $v$, it admits a function $f$ that satisfies every vertex except $v$ and $|\bigcup_{u\in N[v]}{f(u)}| = 4$ by Lemma \ref{exceptional are almost satisfied}.  
Since both graphs can't be $C_4$ (since  $C_4\cdot C_4$ is forbidden and two $C_4$'s joined by a path are prevented by Lemma \ref{extend a tailed cycle}), we are done.

The proof of this lemma is organized as follows.
We first prove structure properties of $G$ in Claims 1-4.
And the rest of the proof is dedicated to a construction of a configuration function of $G$.
It will lead to a contradiction.

\noindent {\bf Claim 1:} $G$ contains no $C_4$'s.

\noindent {\bf Proof of Claim 1:}
Suppose there is a cycle $C=v_1v_2v_3v_4v_1$ of four vertices in $G$.
If there is only one vertex, say $v_1$, in $C$ of degree at least three in $G$, then it is a cut-vertex which is impossible.

Hence there are two vertices in $C$ of degree at least three.
We may assume that the two vertices are $v_1$ and $v_3$.
Let $G'= G \setminus \{v_2\}$.
If $G'$ is configurable, then there is a configuration $f$ on $G'$, and we can extend $f$ to $G$ by assigning $f(v_2) = f(v_4)$, contradicting the assumption that $G$ is not configurable.
Note that $G'$ is a connected graph of maximum degree at most five and of minimum degree at least two with no two vertices of degree at least three adjacent.
Since $\lvert V(G') \rvert + \lvert E(G') \rvert < \lvert V(G) \rvert + \lvert E(G) \rvert$, $G'$ is isomorphic to $C_4$, $C_7$, $C_4\cdot C_4$ or $K_{2,3}$.
If $G'$ is isomorphic to $C_4$, then $G$ is isomorphic to $K_{2,3}$.
If $G'$ is isomorphic to $C_7$, then $G$ is isomorphic to a graph obtained from $C_4$ by adding a path of length five, so $G$ is configurable by Lemma \ref{add a path to cycle}.
If $G'$ is isomorphic to $K_{2,3}$, then $G$ is $K_{2,4}$, and it is configurable by Lemma \ref{small case 2}.
So $G'$ is isomorphic to $C_4\cdot C_4$.
Since $v_4$ is a vertex of degree two and it is a common neighbor of $v_1$ and $v_3$, we have that either $v_1$ or $v_3$ is the vertex of degree four in $C_4\cdot C_4$.
So $G$ can be obtained from adding a path of length four to $K_{2,3}$, so $G$ is configurable by Lemma \ref{add a path to small graph}.
$\Box$

\noindent {\bf Claim 2:} If $P$ is a path whose ends are of degree at least three in $G$ and whose internal vertices are of degree two in $G$, then the number of internal vertices is at most two.

\noindent {\bf Proof of Claim 2:}
If the number of internal vertices of $P$ is at least four, then consider the graph $H$ which is obtained from $G$ by
replacing three consecutive degree two vertices in $P$ by an edge.
If $H$ is configurable, $G$ is also configurable by Lemma \ref{degree two path}.
So $H$ is $C_4$, $C_7$, $C_4\cdot C_4$ or $K_{2,3}$.
But in this case, $G$ can be obtained from $C_4$ by attaching a path of order at least three, so $G$ is configurable by Lemma \ref{add a path to cycle}.
If the number of internal vertices of $P$ is three, then let $H'$ be the graph obtained from $P$ by deleting all internal vertices of $P$.
Again, $G$ is configurable by Lemma \ref{extendable} if $H'$ is configurable.
So $H'$ is $C_4$, $C_7$, $C_4\cdot C_4$ or $K_{2,3}$.
However, $G$ is configurable by Lemma \ref{add a path to cycle} and Lemma \ref{add a path to small graph} in this case.
$\Box$

\noindent {\bf Claim 3:} There are no induced $(\alpha, \beta)$-stars $S$ in $G$, where 
$\alpha + \beta \geq 3$, and $\alpha + 3 \beta \leq 9$ or $(\alpha, \beta) = (1,3)$, such that $G \setminus S$ has minimum degree at least two.

\noindent {\bf Proof of Claim 3:}
Suppose there is an induced $(\alpha, \beta)$-star $S$, where $\alpha + \beta \geq 3$, 
and $\alpha + 3 \beta \leq 9$ or $(\alpha, \beta) = (1,3)$, such that $G \setminus S$ has minimum degree at least two.
Subject to this constraint, assume that $\alpha + \beta$ is as small as possible.
Let $G' = G \setminus S$, and $M_1, M_2 ,..., M_k$ be components of $G'$.
If every component of $G'$ is configurable, then $G$ is configurable by Lemma \ref{extendable}.
So there is a component of $G'$ which is not configurable, and hence this component is isomorphic to $C_4$, $C_7$, $C_4\cdot C_4$ or $K_{2,3}$ by the minimality of $G$.
But $G$ contains no $C_4$'s by Claim 1, so the component is isomorphic to $C_7$.
Without loss of generality, we may assume that $M_1$ is isomorphic to $C_7$ and write $M_1 = v_1v_2 ... v_7 v_1$.

If $M_1$ contains exactly one vertex of degree at least three in $G$, then $G$ is configurable by Lemma \ref{extend a tailed cycle}, a contradiction.
If $M_1$ contains exactly two vertices of degree at least three in $G$, then there is a path of length at least four whose ends are of degree at least three in $G$ and whose internal vertices are of degree two in $G$, contradicting Claim 2.
Hence there are three vertices in $M_1$ of degree at least three in $G$, and we may assume that they are $v_1, v_3, v_5$.
Furthermore, if all $v_1, v_3, v_5$ have degree at least four in $G$, then $\alpha + \beta \geq 6$.
Since $\alpha+3\beta \leq 9$, we have that $\beta \leq 1$ and $G$ contains a $C_4$, contradicting Claim 1.
So at least one of $v_1, v_3, v_5$, say $x$, has degree three in $G$.
Note that there is an $(\alpha, \beta)$-star with center $x$ and $\alpha + \beta =3$ such that the graph 
obtained from $G$ by deleting this $(\alpha, \beta)$-star is still of minimum degree at least two, so $S$ must also have that $\alpha + \beta =3$ by the minimality of $\alpha + \beta$.
So $G'$ is $C_7$ as $\alpha + \beta =3$.
In other words, $(\alpha, \beta) = (0,3), (1,2), (2,1)$ or $(3,0)$.

If $(\alpha, \beta)=(3,0)$, then $G$ can be obtained from $C_6$ by attaching a $(2,1)$-star, so $G$ is configurable by Lemma \ref{extendable}.
So this is not a $(3,0)$-star.
If $(\alpha, \beta)= (0,3)$, then $G$ is configurable since it can be obtained from $C_8$ by attaching a $(1,2)$-star.
If $(\alpha, \beta) = (1,2)$, then $G$ is configurable since $G$ can be obtained from $C_8$ by attaching either a $(2,1)$-star or $(3, 0)$-star.
So $(\alpha, \beta)=(2,1)$.
Let $V(S)=\{a,b,c,d_1, d_2\}$ and $E(S) = \{ab,ac,ad_1, d_1d_2\}$.
If $d_2$ is adjacent to $v_1$ or $v_5$, then $G$ is configurable since it can be obtained from $C_6$ by attaching a $(1,2)$-star.
So $d_2$ is adjacent to $v_3$.
Hence there is a configuration on $G$ defined as $f(v_1) = \{1,2\}, f(v_2) = \{4,5\}, f(v_3) = \{1,3\}, f(v_4) = \{4,5\}, f(v_5) = \{1,2\}, f(v_6) = \{3,4\}, f(v_7) = \{3,5\}, f(a) = \{1,3\}, f(b) = f(c) = \{4,5\}, f(d_1) = \{2,5\}, f(d_2) = \{2,4\}$.
This proves Claim 3.
$\Box$

\noindent {\bf Claim 4:} $G$ contains no $C_6$ with exactly two vertices of degree at least three that are diagonally opposite on the cycle.

\noindent {\bf Proof of Claim 4:}
Let $C = v_1 v_2 ... v_6 v_1$ be a cycle of order six with $v_1$ and $v_4$ the two vertices of degree at least three in $G$.
Since $G$ has no adjacent vertices whose degrees are at least three, $v_5$ and $v_6$ have degree two in $G$.
Let $G'$ be the graph obtained by deleting $v_5, v_6$ from $G$, so $G'$ is a graph of minimum degree at least two, maximum degree at most five, and there are no adjacent vertices whose degrees are at least three.
If $G'$ is not configurable, then $G'$ is $C_4$, $C_7$, $C_4\cdot C_4$ or $K_{2,3}$ by the minimality of $G$.
However, $G$ contains no $C_4$'s, so $G'$ is $C_7$ and it contains at most two vertices whose degrees in $G$
 are at least three.
Hence, there is a path of order at least five whose internal vertices are of degree two, which contradicts to Claim 2.
Consequently, $G'$ is configurable and there is a configuration $f$ on $G'$, and we can extend $f$ to $G$ by defining $f(v_5) = f(v_3)$ and $f(v_6) = f(v_2)$.
$\Box$

We now construct a configuration on $G$.
Construct a graph $H$ as follows: the vertices of $H$ are the vertices of degree at least three in $G$, 
and $xy$ is an edge in $H$ if $x$ and $y$ have a common neighbor in $G$.

\noindent {\bf Claim 5:} The maximum degree of $H$ is at most two.

\noindent {\bf Proof of Claim 5:}
Suppose there is a vertex $x$ of degree at least three in $H$.
Let $x_1, x_2, ..., x_k$ be the vertices of degree at least three such that there exist $x$-$x_i$ paths of length two or three.
Then the internal vertices of those $x$-$x_i$ paths together with $x$ form an $(\alpha, \beta)$-star $S$ with $\alpha \geq 3$.
On the other hand, $\alpha + \beta$ is at most five since $G$ is of maximum degree at most five.
So $S$ is an $(\alpha, \beta)$-star with $\alpha + 3 \beta \leq 9$.
By Claim 3, $G \setminus S$ is not of minimum degree at least two.
So the degree of $x_i$ in $G \setminus S$ is at most one, for some $i=1,2,..., k$.
Since $G$ contains no $C_4$'s and $C_6$ with exactly two diagonal vertices of degree at least three in $G$, the degree of $x_i$ is exactly three.
So there is an $(\alpha', \beta')$-star $S'$ centered at $x_i$ with $\alpha' + 3 \beta' \leq 9$ such that $G \setminus S'$ is of minimum degree two since $\alpha \geq 3$, which contradicts Claim 3.
Hence, the maximum degree of $H$ is at most two.
$\Box$

By Claim 5, $H$ is a disjoint union of isolated vertices, paths and cycles.
Let $H^2$ be the graph obtained by adding edges $xy$ to $H$ for each pair of two vertices $x,y$ which have distance exactly two between them in $H$, and then deleting multiple edges and loops.
So $H^2$ has maximum degree at most four.
Let $H'$ be the graph that is obtained by deleting an edge which is in $H^2$ but not in $H$ from each component of $H^2$ isomorphic to $K_5$.
Hence, $H'$ is $4$-colorable by Brooks' Theorem.
Let $c: V(H') \rightarrow \{1,2,3,4\}$ be a proper $4$-coloring of $H'$ such that $c(v)=1$ for each isolated vertex $v$ in $H$.
Note that $H^2$ contains a component which is isomorphic to $K_5$ if and only if the component in $H$ is isomorphic to $C_5$.

Define a function $f: V(H) \rightarrow [5]^2$ as $f(v) = \{c(v), 5\}$ for every vertex $v$ in $H$.
Let $U$ be the set of vertices $u$ such that $u$ is a common neighbor of two vertices of degree at least three in $G$.
Since no two vertices of degree at least three are adjacent, every vertex in $U$ is of degree two in $G$.
Now, we shall extend $f$ to $V(H) \cup U$ by defining $f(u) = \{1,2,3,4,5\} \setminus (f(x) \cup f(y))$ for each vertex $u$ in $U$, where $x,y$ are the two neighbors of $u$ in $G$.
Note that if $x,y$ are the two neighbors of a vertex $u$ in $U$, then $c(x) \neq c(y)$ since $H'$ contains all edges in $H$, so $\lvert f(x) \cup f(y) \rvert =3$, and $f$ is well-defined on $V(H) \cup U$.
It is clear that $\bigcup_{w \in N[u]} f(u) = \{1,2,3,4,5\}$ for each $u \in U$.
Furthermore, if $v$ is a vertex with degree at least two in $H$, and $v$ is not in a component of $H$ isomorphic to $C_5$, then neighbors of $v$ in $H$ receive different colors under $c$, so $u$ is satisfied.
Similarly, for each component of $H$ which is isomorphic to $C_5$, there is a vertex $w$ such that $\lvert \bigcup_{u \in N[w] \cap (V(H) \cup U)} f(u) \rvert = 4$ and each other vertex is satisfied.

Let $W$ be the set of vertices $w$ that are not satisfied.
So each vertex in $W$ is either an isolated vertex in $H$, an end of a maximal path in $H$, or a vertex in a component of $H$ which is isomorphic to $C_5$.
Let $X = \{w \in W: w$ is an isolated vertex in $H\}$, and let $Y$ be the set $W \setminus X$.
Notice that $\lvert \bigcup_{u \in N[w] \cap (V(H) \cup U)} f(u) \rvert =4$ when $w$ is in $Y$.
Now, construct a graph $L$, where $V(L)$ is equal to $V(H)$, and two vertices $x,y$ in $L$ are adjacent if there is a $x$-$y$ path of length three in $G$.
Note that since no vertices of degree at least three are adjacent, the internal vertices of every $x$-$y$ path of length three in $G$ are of degree two for each $xy \in E(L)$.

\noindent {\bf Claim 6:}
If $w$ is in $X$, then the degree of $w$ in $L$ is at least four.
If $w$ is in $Y$, then the degree of $w$ in $L$ is at least two.

\noindent {\bf Proof of Claim 6:}
Let $w$ be a vertex in $X \cup Y$.
Let $x_1, x_2 ,..., x_k$ be vertices of degree at least three in $G$ such that there are $w$-$x_i$ paths in $G$ of length two or three for each $i=1,2,..., k$.
Then the internal vertices of those $w$-$x_i$ paths together with $w$ form an $(\alpha, \beta)$-star $S$.

Suppose $w \in X$.  Then $\alpha = 0$ and there is at most one path between $w$ and each $x_i$ since otherwise we violate Claim 4.  But then $G\backslash S$ has minimum degree two, so by Claim 3, $\beta \geq 4$, so the degree of $w$ in $L$ is at least four.

Suppose $w \in Y$ and that $\beta \leq 1$.  If $w$ was not in a $C_5$ in $H$, then $\alpha = 1$, so the degree of $w$ is only two.  So we must have that $w$ was in a $C_5$ in $H$, so $\alpha = 2$.  Removing $S$ must create a vertex of degree one by Claim 3, say $x_1$.  
So $x_1$ must have degree three and be part of a 5-cycle $D$ in $G$ with $w$.  
Since $w$ is in a $C_5$ in $H$, $G$ must have that $x_1$ has a path of length two to another vertex of degree at least three in $G$ and that the graph $H'$ obtained from $G$ by removing $D$ and the two degree two vertices that are adjacent to vertices of $D$ is connected and of minimum degree two.  
If $H'$ is configurable, then by Lemma \ref{add a C_5 with two tails}, $G$ would be as well, so $H'$ must be $C_7$ which is impossible since it has at least one degree three vertex since $G$ has at least five degree three vertices since $w$ is in a $C_5$ in $H$.
$\Box$

By Lemma \ref{orientation}, $L$ then has an orientation in which each vertex of $X$ has in-degree at least two and every vertex in $Y$ has in-degree at least one.  We use this to extend $f$ to satisfy every vertex in $G$.  Each edge in $L$ corresponds to a path of length three, $x, v_1, v_2, y$ in $G$ (where $x$ is the tail of the edge in $L$).  For each of these paths, let $a, b$ be two colors not in $\bigcup_{u\in N(x)}{f(u)}$ (if that many colors exist, otherwise arbitrarily add colors not in $f(x)$).  Then assign $f(v_1) = (a, b)$ and $f(v_2)$ as given by Lemma \ref{P_2 helper}.

Clearly at the end of this process each vertex of degree two is satisfied.  Each vertex not in $X$ or $Y$ was already satisfied.  Each vertex in $X$ was the tail of two edges in $L$, so sees up to four new colors, so is certainly satisfied.  Each vertex in $Y$ was only missing at most two colors, but was the tail of at least one edge in $L$, so is now satisfied.
\end{pf}

\section{Main theorem}

We now prove Theorem~\ref{thm:main}, which we restate in equivalent form.

\begin{theorem}
If $G$ is a connected graph of minimum degree at least two with no induced subgraph isomorphic to $K_{1,6}$, and $G$ is not isomorphic to a member of $\{C_4,C_7,C_4\cdot C_4,K_{2,3},G_i: 1 \leq i \leq 4\}$, then $G$ is configurable.
\end{theorem}

\begin{pf}
We first prove the theorem for graphs on at most six vertices.
It is easy to see that the theorem holds if $\lvert V(G) \rvert \leq 4$, so we assume that $5 \leq \lvert V(G) \rvert \leq 6$.
If $G$ contains $C_6$, then $C_6$ is a spanning subgraph of $G$.
Since $C_6$ is configurable, $G$ is configurable.
So we may assume that $G$ does not contain $C_6$.
If $G$ contains $C_5$, then $G$ contains a spanning subgraph that is obtained from $C_5$ by attaching a path on one vertex.
Since $G$ does not contain $C_6$, $G$ is configurable by Lemma \ref{5 6 cycle}.
Hence, we may assume that the longest cycle in $G$ has length at most four.

Assume that $G$ contains $C_4$.
Since $\lvert V(G) \rvert \leq 6$, $G$ is $2$-edge-connected.
So $G$ contains a spanning subgraph that can be obtained from $C_4$ by consecutively attaching paths.
If the first path we attached contains two vertices, then since $G$ has no cycle of length greater than four, $G$ contains a spanning subgraph that can be obtained from a triangle by attaching a path on three vertices and hence is configurable by Lemma \ref{add a path to cycle}.
If the first path we attached has only one vertex, then since $G$ does not contain $C_5$, $G$ contains a spanning subgraph that can be obtained from a triangle by attaching two paths on one vertex to different vertices or from $K_{2,3}$ by attaching a path on one vertex, so we are done by Lemmas \ref{small extendable} and \ref{add a path to small graph}.

Therefore, we may assume that every cycle in $G$ is a triangle.
If $G$ is $2$-edge-connected, then $G$ can be obtained from $C_3$ by attaching a path on two vertices and hence is configurable by Lemma \ref{small extendable}.
If $G$ is not $2$-edge-connected, then $G$ contains two disjoint triangles as a spanning subgraph, and hence $G$ is configurable.
This proves that the theorem holds for graphs on at most six vertices.
	
We now proceed by induction on $\lvert V(G) \rvert + \lvert E(G) \rvert$.
We have shown the theorem holds for graphs on at most six vertices, so we may assume that the order of $G$ is at least seven.

Suppose there is a vertex $v$ of degree two in $G$ such that $v$ is in a $C_4=vabcv$ with degree of $b$ also two.  
Note that $G_i$ contains a spanning cycle of length seven for $1 \leq i \leq 4$.
Suppose that the degree of $a$ is also two.
If $c$ is not of degree three, then $G$ is obtained by attaching a path on three vertices to a configurable graph or an exceptional graph, so $G$ is configurable by Lemmas \ref{extendable}, \ref{add a path to cycle} and \ref{add a path to small graph}.  
If $c$ is of degree three, then $G$ is obtained from a $C_4$ and a graph by attaching a path, where the ends of the path are adjacent to vertices in different components.
Then $G$ is configurable by Lemmas \ref{extend a tailed cycle} and \ref{connectivity}.
So we may assume that $a$ and $c$ have degree at least three.

So we consider $G \setminus v$.  If it has a configuration $f$, then $G$ is configurable since we may extend $f$ to $V(G)$ by assigning $f(v)=f(b)$.
As the order of $G$ is at least seven, $G \setminus v$ is not configurable only if $G \setminus v$ is $C_4\cdot C_4$ or contains a spanning cycle of length seven.
However, it is not hard to see that if $G \setminus v$ is $C_4\cdot C_4$ or contains a spanning cycle of length seven, then $G$ contains a spanning subgraph that can be obtained either from $C_4 \cdot C_4$ by attaching a path on one vertex or from $C_4$ by attaching a path on four vertices, so $G$ is configurable by Lemmas \ref{add a path to cycle} and \ref{add a path to small graph}.
Hence we may assume that no four cycle has two vertices of degree two opposite one another.

Suppose there were three vertices $x,y,z$ in $G$ such that $x,y,z$ form a triangle in $G$ and the degree of $y$ and $z$ were exactly two.
Assume that $x$ is not of degree three.
By the induction hypothesis, Lemma \ref{small extendable} and Lemma \ref{add a path to small graph}, $G$ is configurable if $G \setminus \{y,z\}$ is not $C_4$ or contains $C_7$ as a spanning subgraph. 
But if $G \setminus \{y,z\}$ is $C_4$ or contains $C_7$ as a spanning subgraph, then $G$ contains a spanning subgraph that can be obtained from $C_3$ by attaching a path with order at least three, so $G$ is still configurable by Lemma \ref{add a path to cycle}.
Similarly, if $x$ is of degree three, then $G$ is configurable by Lemma \ref{connectivity}.
Hence, we may assume that $G$ has no triangles with two vertices of degree three.

Let $G'$ be a spanning subgraph of $G$ such that the minimum degree of $G'$ is at least two and $G'$ satisfies the following:
\begin{enumerate}
	\item $\lvert E(G') \rvert$  is as small as possible, 
	\item subject to that, the number of triangles in $G'$ is as small as possible, and
	\item subject to that, the number of components in $G'$ which are isomorphic to $C_4\cdot C_4$ or $K_{2,3}$ is as small as possible.
\end{enumerate}
We shall prove the following claims.

Note that by the minimality of $E(G')$, there are no two vertices of degree at least three adjacent to one another.

\noindent {\bf Claim 1:} The maximum degree of $G'$ is at most five.

\noindent {\bf Proof of Claim 1:}
Suppose that there is a vertex $v$ of degree at least six in $G'$.
As $G$ is $K_{1,6}$-free, there are two vertices $x,y$ adjacent to $v$ in $G'$ with $x$ adjacent to $y$ in $G$.
Since the degree of $v$ is at least three, $x$ and $y$ must have degree two in $G'$.
If $xy \not \in E(G')$, then the graph obtained by deleting $xv, yv$ from $G'$ and then adding $xy$ into $G'$ is still a spanning subgraph of $G$ with minimum degree at least two, but it has fewer edges.
So $xy \in E(G')$, in other words, $v,x,y$ form a triangle in $G'$.
Since $x,y,v$ form a triangle in $G$ and the degree of $v$ is at least three, at least one of $x$ and $y$ has degree at least three in $G$.
We may assume that the degree of $x$ in $G$ is at least three, and $u$ is a neighbor of $x$ in $G$ other than $y$ and $v$.
As $xy, vx \in E(G')$ and the degree of $x$ is two in $G'$, $xu \not \in E(G')$.
So the graph obtained by deleting $xv$ and adding $xu$ has the same number of edges but it has fewer triangles than $G'$, a contradiction.
$\Box$


%

Since every component of $G'$ is a connected graph of minimum degree at least two and of maximum degree at most five, and no vertices of degree at least three in $G'$ are adjacent to one another, every component of $G'$ is configurable except those that are isomorphic to $C_4$, $C_7$, $C_4\cdot C_4$, or $K_{2,3}$ by Lemma \ref{min degree 2 max degree 5}.
Also, it follows by a simple case checking that if a graph not containing $C_7$ as a spanning subgraph contains $C_4$, $C_4\cdot C_4$ or $K_{2,3}$ as a spanning subgraph but not as an induced subgraph, then it is also configurable.


Now, we show that $G$ is configurable.
If $\lvert V(G) \rvert=7$ but $G$ is not configurable, then $G$ contains $C_7$ as a spanning subgraph.
We denote the $C_7$ by $v_0v_2...v_6$.
If there exists $i$ with $0 \leq i \leq 6$ such that $v_iv_{i+2}$ is an edge, where the index is computed modulo seven, then $G$ contains a spanning subgraph that can be obtained from $C_3$ by adding a path on four vertices, so $G$ is configurable by Lemma \ref{add a path to cycle}.
Since $G$ is not $C_7$ or $G_1$, $G$ contains at least nine edges.
If there exists $i$ with $0 \leq i \leq 6$ such that $v_iv_{i+3}$ and $v_{i+1}v_{i+5}$ are edges of $G$, then $G$ is configurable by Lemma \ref{5 6 cycle}.
So $G$ contains $G_2$ or $G_3$ as a subgraph but not an induced subgraph since $G$ is not $G_2$ or $G_3$.
In addition, adding an edge to $G_2$ or $G_3$ makes it configurable unless it creates $G_4$.
But adding an edge to $G_4$ makes it configurable.
This proves that $G$ is configurable if $G$ contains at most seven vertices.
So we may assume that $G$ has at least eight vertices.

Let $H$ be a maximal configurable subgraph of $G$ induced by a union of components of $G'$.  
Suppose that $H$ is empty.
Since $G$ contains at least eight vertices, $G'$ contains at least two components.
Let $H_1,H_2$ be two components of $G'$ adjacent in $G$ and $v_i$ be a vertex of $H_i$ adjacent in $G$ to $H_{3-i}$ for $i=1,2$.
By Lemma \ref{exceptional are almost satisfied}, for each $C_4,C_7,C_4 \cdot C_4$ and $K_{2,3}$, and for each of its vertices $v$, there exists a function $f$ mapping the vertices to $[5]^2$ satisfying every vertex except possibly $v$, and $v$ is missing at most two colors.
Let $f_1,f_2$ be such a function defined on $V(H_1)$ and $V(H_2)$, respectively, such that $v_1,v_2$ are the only vertices missing some colors.
Therefore, we can permute the colors in $f_1$ and $f_2$ such that $f_i(v_i)$ contains the colors which $v_{3-i}$ missed for $i=1,2$.
This proves that the subgraph of $G$ induced by $V(H_1) \cup V(H_2)$ is configurable, so $H$ is not empty.

If $H \neq G$, then let $C$ be a component of $G'$ disjoint from $H$ but adjacent in $G$ to $H$.  
By Lemma \ref{exceptional are almost satisfied}, for every $v \in V(C)$, there exists a function $f$ mapping the vertices to $[5]^2$ satisfying every vertex except possibly $v$, and $v$ is missing at most two colors.
Therefore, the subgraph of $G$ induced by $V(H) \cup V(C)$ is configurable by Lemma \ref{connectivity}, contradicting the maximality of $H$.
This proves that $H=G$ and $G$ is configurable.
\end{pf}

\bigskip

Note that our proof gives a polynomial-time algorithm to find a configuration of an $n$-vertex graph $G$ 
if $G$ is a $K_{1,6}$-free graph 
of minimum degree at least two, and no component of $G$ is isomorphic to  $C_4$, $C_7$, $C_4\cdot C_4$ or $K_{2,3}$.

Now we shall show that the hypothesis that $G$ be $K_{1,6}$-free cannot be replaced by assuming that
$G$ be $K_{1,9}$-free. We do so by exhibiting infinitely many examples that contain no induced $K_{1,9}$ \
but are not configurable.
Let $H'$ be the graph obtained from $K_5$ by replacing each edge $xy$ by two internally disjoint paths 
$xu_{xy}y$ and $x v_x v_y y$, and $H$ be the graph obtained from $H'$ by deleting $v_a$ and $v_b$, 
where $a$ and $b$ are two distinct vertices in the original $K_5$.
So the maximum degree of $H$ is eight, and there are exactly two vertices which have degree seven.
Suppose that $H$ is configurable and $f$ is a configuration on $H$.
If $x,y$ are distinct vertices in the original $K_5$, then $f(x) \neq f(y)$ for otherwise $\bigcup_{z \in N[u_{xy}]} f(z) \neq \{1,2,3,4,5\}$, and $f(x) \cap f(y)$ is nonempty for otherwise $\bigcup_{z \in N[v_x]} f(z)$ or $\bigcup_{z \in N[v_y]} f(z)$ is not $\{1,2,3,4,5\}$.
But if $S$ is a subset of $[5]^2$ such that every two members of $S$ have a nonempty intersection, then the size of $S$ is at most four, so $f(a) = f(b)$.
However, this implies $\bigcup_{w \in N[u_{ab}]} f(w) \neq \{1,2,3,4,5\}$, a contradiction.
Hence, $H$ is not configurable.
For any positive integer $k$, let $H_1, H_2, ..., H_k$ be graphs, where each of them is isomorphic to $H$, and $a_i,b_i$ are the two vertices of degree seven of $H_i$ for each $i=1,2,..., k$.
Let $G$ be the graph obtained from $H_1\cup H_2\cup\cdots\cup H_k$ by  adding the edges $b_i a_{i+1}$ for all 
$i=1,2,...,k-1$ and $b_k a_1$, so $G$ is of maximum degree eight but not configurable.

\begin{figure*}[ht]
\begin{center}
\includegraphics[scale=0.7]{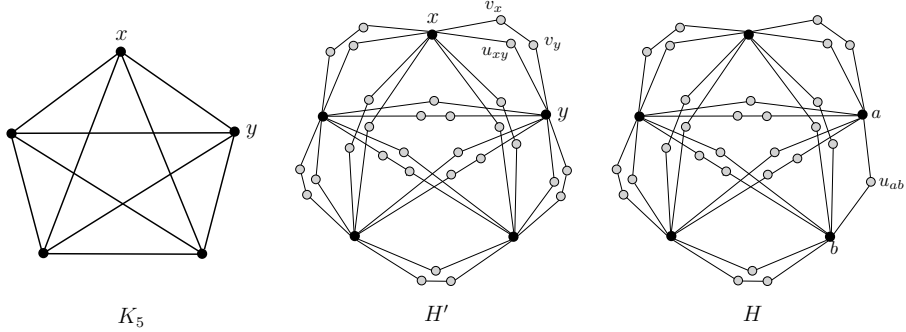}
\caption{A complete graph $K_5$. $H'$ is obtained by replacing every edge $xy \in E(K_5)$, by disjoint paths $xu_{xy}y$ and $x v_x v_y y$. $H$ is obtained from $H'$, by deleting $v_a$ and $v_b$, from two distinct vertices $a$ and $b$. Note that every vertex of $H$ 
that belongs to the original $K_5$ has degree eight, except $a$ and $b$, that have degree seven.}
\label{fig:K19_free}
\end{center}
\end{figure*}

On the other hand, one might ask whether we can get rid of the assumption about forbidden subgraphs by assuming the minimum degree is large.
However, the following examples show that for every integer $k>0$, there is a graph $G$ with minimum degree $k$
that is not configurable.
Let $n=10k-9$, let $B$ be a set of size $n$, and let $A$ be the set of all $k$-element subsets of $B$.
Let $G$ be the graph with vertex-set $A\cup B$ in which a vertex $S\in A$ is adjacent to each of its elements.
By the pigeon hole principle there is a set $S$ in $A$ such that $f(b)$ are the same for all $b \in S$. 
But this implies that $\lvert \bigcup_{v \in N[S]}f(v) \rvert \leq 4$, a contradiction.
So $G$ is not configurable.

\bigskip

\noindent {\bf Acknowledgement.} The authors thank Wayne Goddard for pointing out that the list of the exceptional graphs of the main theorem was incomplete in the first version of this paper.

\baselineskip 11pt
\vfill
\noindent
This material is based upon work supported by the National Science Foundation.
Any opinions, findings, and conclusions or
recommendations expressed in this material are those of the authors and do
not necessarily reflect the views of the National Science Foundation.

\end{document}